\documentclass{article}
\usepackage{graphicx}
 \usepackage{mathptmx}
\usepackage{amsmath, amstext,amssymb,amsfonts}
\usepackage[english]{babel}
\usepackage{delarray}
\usepackage{mathptmx}
\usepackage{amsmath, amstext,amssymb,amsfonts}
\usepackage[english]{babel}
\usepackage{delarray}

\DeclareFontFamily{U}{mathx}{\hyphenchar\font45}
\DeclareFontShape{U}{mathx}{m}{n}{
      <5> <6> <7> <8> <9> <10>
      <10.95> <12> <14.4> <17.28> <20.74> <24.88>
      mathx10
      }{}
\DeclareSymbolFont{mathx}{U}{mathx}{m}{n}
\DeclareFontSubstitution{U}{mathx}{m}{n}
\DeclareMathAccent{\widecheck}{0}{mathx}{"71}
\DeclareMathAccent{\wideparen}{0}{mathx}{"75}

\newcommand{\R}{\mathbb R}      \newcommand{\Co}{\mathbb C}  \newcommand{\usi}{U_{\sigma}}   \newcommand{\Z}{\mathbb Z}
       \newcommand{\B}{B(\R)}   \newcommand{\Le}{L^1({\mathbb R})}  \newcommand{\M}{M({\mathbb R})}
   \newcommand{\us}{U_{\sigma}}  \newcommand{\bp}{B_{\sigma}^{p}}   \newcommand{\be}{B_{\sigma}^{1}}

\newtheorem{theorem}{Theorem}[section]
\newtheorem{lemma}[theorem]{Lemma}
\newtheorem{prop}[theorem]{Proposition}
\newtheorem{cor}[theorem]{Corollary}
\newtheorem{con}[theorem]{Conjecture}
\newtheorem{remark}[theorem]{Remark}

\numberwithin{equation}{section}

\date{}
\title{A note on uniqueness of extension for characteristic functions}
\author{\large{ Saulius  Norvidas}}
\date{\footnotesize Institute of Data Science and Digital Technologies, Vilnius University, \\ Akademijos str. 4, Vilnius LT-04812, Lithuania\\
 ({\rm{e-mail: norvidas{@}gmail.com}})}
\begin{document}

\maketitle
 {{ {\bf Abstract}}}
  Let $f:\R\to\Co$ be the characteristic function of a probability measure. We study the following question: Is it true that for any closed interval $I$ on  $\R$,
which does not contain the origin, there exists a characteristic function $g$ such that $g$ coincides with $f$   on $I$  but $g \not\equiv f$  on $\R$?



\section{ Introduction }

{\large{
{\large{
Let  $M(\R)$   be  the set of bounded complex-valued regular Borel measures on the real line $ \R$.   For  $\mu \in  M(\R)$,  we define its Fourier transform  by
\[
\hat{\mu}(t)=\int_{\R}e^{-it x} d\mu(x), \quad t\in\R.
\]
If  $\mu\in\M$  is   non-negative and $\|\mu\|=1$, then in the language of probability theory,   $\mu$ and    $f(t)= \hat{\mu}(t)$  are called a probability measure and its   characteristic function, respectively. Any characteristic function $f$ satisfies $f(0)=1$ and
\begin{equation}
f(-t)=\overline{f(t)}, \quad t\in\R.
\end{equation}
  Let $\Le\subset \M$ be the usual  Lebesgue space  of  absolutely continuous measures (with respect to  Lebesgue measure $dx$). If $\varphi\in L^1(\R)$, $\|\varphi\|_{L^1(\R)}=1$, and $\varphi\ge 0$  on $\R$,  then $\varphi$  is called a probability density function, or (probability) density for short. Note that  if $\varphi$ is a measurable function,  then we write $\varphi\ge 0$ on $\R$, if $\varphi\ge 0$   almost everywhere, i.e.,  $dx$-almost everywhere on $\R$.

A  complex-valued function $f$ defined on $(-A, A)$, $0<A\le\infty$,  is said to be   positive definite if
\[
\sum_{j, k=1}^{n} f(t_j-t_k)c_j\overline{c}_k\ge 0
\]
 for  all finite sets  $c_1,\dots, c_n\in\Co$ and  $t_1, \dots, t_n \in (-A/2, A/2)$.  The Bochner theorem (see, e.g., \cite[p. 71]{11}) states that  a continuous  function $f : \R\to \Co$ is characteristic  function if and only if $f$  is   positive definite on $\R$ and $f(0)=1$.

In this paper, we shall be concerned with a  problem of uniqueness for extensions of characteristic  functions. By a classical results  of Krein \cite{8} (see also  Akutowicz \cite{1}, Kahane \cite{7}, Levin \cite{10}, and Sasv\'{a}ri \cite{14}),    for any neighbourhood   $U$  of the origin and a  continuous positive definite function $f$ on $U$,  there  exists a continuous  extension of $f$  to  $\R$ so as to retain the positive definite property. In general such an extension is not unique  (see, e.g., \cite[p. 89]{11}).   A question concerning  the uniqueness of that  extensions  was also studied. Conditions for the uniqueness can be found in \cite{8}.

We consider a certain  extension problem in the case if $U$  is a neighbourhood of infinity. Our study is  initiated  by the   question  posed by  N.G. Ushakov (\cite[p. 276]{16})  :  Is it true that for any  interval $[a, b]\subset \mathbb{R} $, $0\notin [a,b]$,  there exists the characteristic function $f$  such that  $f \not\equiv e^{-t^2/2}$, but  $ f (t) = e^{-t^2/2}$ for all $t\in [a, b]$?
 A positive answer to this question  was given by T. Gneiting in \cite{4}. It follows from the following result:
\begin{theorem}\cite[p.  360]{4}.
 Let $f : \R\to\Co$  be the characteristic function of a  distribution with a continuous and strictly positive density.  Then there exists, for each $a>0$,   a characteristic function $g$   such that $f(t)=g(t)$ if $t=0$ or $|t|\ge a$ and  $f(t)\neq g(t)$ otherwise.
 \end{theorem}
Repeating the author's definition \cite[p.  361]{4}, we say that a characteristic function $f$ has the substitution property if, for any $a>0$, there exists a characteristic function $g$ such that $f(t)=g(t)$ for $|t|> a$ but   $f\not\equiv g$. Also in \cite[p.  361]{4} we find the following: {\it One might conjecture that any characteristic function with an absolutely continuous component has the substitution property}. This conjecture  fails in general.   Indeed, for  $a>0$,
\[
\varphi_a(x)  =\begin{array}\{{rl}.
\frac{2(a- 2|x|)}{a^2}, &  |x|\le \frac{a}2, \\
0, & \mbox{otherwise,}
\end{array}
\]
is the density of the triangular probability distribution. Let $\sigma>0$ and assume  that $g$ is a characteristic function such that $g = \widehat{\varphi_a}$ on the neighbourhood of infinity  $\R\setminus [-\sigma, \sigma]$. If $a$ satisfies $a\le \pi/\sigma$, then  $g= \widehat{\varphi_a}$ on the whole $\R$ (see \cite[Example 1, p. 238]{12}). Below we improve this estimate in our  Corollary 1.4.  Moreover, in  Theorem 1.5 and Proposition 1.8 we show that the exact estimate is $a\le 2\pi/\sigma$.

Let us recall the  notation that will be used throughout this paper.  In view of (1.1), it is enough to study the extensions of characteristic functions only from symmetric neighborhoods of infinity.  Therefore, in the sequel,
\[
U_{\sigma}=\{x\in\R: \ |x|>\sigma\},
\]
 $\sigma>0$, denotes such a neighborhood. Given  a measurable function $\varphi : \R \to \R$, we denote by $S_{\varphi}$ the essential support of $\varphi$. By definition, a point $x\in\R$ belongs to $S_{\varphi}$ if, for any $\varepsilon>0$, the set
\[
(x-\varepsilon, x+\varepsilon)\cap \{t\in\R: \ |\varphi(t)|>0\}
\]
has positive Lebesgue measure.  We call   $N_{\varphi}=\R\setminus S_{\varphi}$  the essential zero set of $\varphi$. The Lebesque measure of  a measurable   $E\subset\R$  will be denoted  by $|E|$.   Let $A$ be  a proper subset   of $\R$.  For functions  $f$ and  $g$ defined on $A$,  we  write $f=g$  on $A$  if $f(x)=g(x)$ for all $x\in A$. If  $f(x)=g(x)$ for all $x\in \R$, then we  write $f\equiv g$.  Finally, let   $\Z$  be the group of integers as usual.

We now formulate our  results.  We start by proving  that it is sufficient to study a smaller class of characteristic functions.
 \begin{prop}
Suppose that $ \mu$ and $\eta$ are  probability measures on $\R$. Let $\mu=\mu_a+\mu_s$ and $\eta=\eta_a+\eta_s$ be the usual Lebesgue decompositions of $ \mu$ and $\eta$ into their  absolutely continuous and singular parts, respectively. If there exists  a $\us$ such that $\widehat{\mu}=\widehat{\eta}$  on $\us$, then $\mu_s=\eta_s$.
\end{prop}
For this reason, we can restrict our  extension problem  to a smaller class  of characteristic functions of probability densities.
The following uniqueness theorem  is the main result  of this paper.
\begin{theorem}
Let  $f : \R\to\Co$  be the characteristic function of a probability   density $\varphi$.    Suppose that there exist  $a>0$  and measurable set $E\subset [0,a)$  such that  $|E|>0$ and
  \begin{equation}
 \Bigl|(E+a\Z)\cap S_{\varphi}\Bigr|=0.
\end{equation}
   Let $g: \R\to\Co$  be a characteristic function,  and let $f=g$ on $\usi$. If
   \begin{equation}
 a \le \frac{2\pi}{\sigma},
\end{equation}
then $f\equiv g$.
\end{theorem}

\begin{cor}
Let $f$, $g$,  and $\varphi$ be the same as that in Theorem 1.3. If $f=g$ on $\usi$ and
\begin{equation}
 \bigl| S_{\varphi}\bigr|< \frac{2\pi}{\sigma},
\end{equation}
then  $f\equiv g$.
\end{cor}

 The statements of Theorem 1.3 and Corollary 1.4   are  sharp in the sense that the right-hand sides of (1.3) and (1.4) cannot be replaced by $2\pi/\sigma+\varepsilon$ for any positive $\varepsilon$. Indeed, we can  show the following.

\begin{theorem}
Suppose that  $f : \R\to\Co$ is the characteristic function of a continuous probability  density $\varphi$. Let $\alpha, \beta\in\R$, $\alpha<\beta$, be such that
\begin{equation}
(\alpha,\beta)\subset  S_{\varphi}.
\end{equation}
 If
\begin{equation}
 \beta-\alpha > \frac{2\pi}{\sigma},
\end{equation}
 then there is a characteristic function $g$  such that $f=g$ on $U_{\sigma}$ but $f\not\equiv g$.
\end{theorem}
It is easy to see that Theorem 1.5  provides the following sufficient conditions under which  a characteristic function has the substitution property.
\begin{cor}
Let  $f : \R\to\Co$ be the characteristic function of a continuous probability  density $\varphi$. If, for any $\delta >0$,  there exists $A=A(\delta)\in \R$ such that
\begin{equation}
(A, A+\delta)\subset S_{\varphi},
\end{equation}
then $f$ has the substitution property.
\end{cor}
\begin{con}\label{6}
Suppose that  $f$ is the characteristic function of a probability measure with a nontrivial absolutely continuous component.  If  $\varphi$ is the density function of this component, then $f$ has the substitution property if and only if  $N_{\varphi}$ does not contain any  lattice $\tau+a\Z$,  $\tau\in\R$, $a>0$.
\end{con}

Of course, if $\varphi$ is continuous, then the sufficiency part of this conjecture follows from Corollary 1.6.

Finally, we show that for the limit case $\beta-\alpha = 2\pi/\sigma$ in (1.4) and  (1.6), there is no an exact answer to the question for  uniqueness of characteristic extension from $U_{\sigma}$.
\begin{prop}
Given $\sigma>0$ and $\alpha\in\R$, let  $\varphi$ be  any  probability density with
\[
S_{\varphi}=\Bigl(\alpha, \alpha+\frac{2\pi}{\sigma}\Bigr).
\]
(i) There exists a $\varphi$ such that for any characteristic function $g:\R\to\Co$ with $\widehat{\varphi}=g$ on $U_{\sigma}$, it follows that $\widehat{\varphi}\equiv g$.

(ii) There exist a $\varphi$ and  the characteristic function $g:\R\to\Co$ such that  $\widehat{\varphi}=g$ on $U_{\sigma}$ but $\widehat{\varphi}\not\equiv g$.
\end{prop}
We conclude this section by presenting our previous paper \cite{12}, where a similar extension problem was studied in the case of continuous density functions. The main result of \cite{12} states that if $\varphi$ is a continuous probability density such that there exist lattices $\Lambda_j=\tau_j+\alpha_j\Z$, $\tau_j\in\R$, $\alpha_j>0$, $\alpha_j\sigma\le 2\pi$, $j=1,2$, $\Lambda_1\cap \Lambda_2=\emptyset$,  and $\varphi$ vanishes on  $\Lambda_1\cup \Lambda_2$, then, for any characteristic function $g:\R\to\Co$ with $g=\widehat{\varphi}$ on  $U_{\sigma}$, we have that $g\equiv\widehat{\varphi}$. It is easy to see that for continuous density $\varphi$ this statement is more general than our Theorem 1.3. On the other hand, the formulation of this statement (as also as its proof) uses substantially the fact that $\varphi$ is continuos.

\section{ PRELIMINARIES AND PROOFS}

 Let $\B=\{\widehat{\mu}: \mu\in\M\}$ denote the Fourier-Stieltjes algebra with the usual pointwise multiplication.
The norm in  $\B$   is inherited  from  $\M$ in such a way
\[
\|\hat{\mu}\|_{\B}:=\|\mu\|_{\M}.
 \]
 We normalize the inverse Fourier transform
\[
\check{\omega}(x)=\frac1{2\pi}\int_{\R}e^{ix t}\omega(t)\,dt
\]
so that the inversion formula $\widehat{(\check{\omega})}=\omega$  is true  for suitable ƒ$\omega\in \Le$.

As usual, we write $S(\R)$ for the  Schwartz space of test functions on $\R$ and $S'(\R)$ for the dual space of tempered distributions.
Let  ƒ$\Omega$ be a closed subset of $\R$. A function ƒ$\omega\in L^p(\R)$,  $1\le p\le\infty$,   is called bandlimited to ƒ$\Omega$  if ƒ$\widehat{\omega}$  vanishes outside $\Omega$. If  $2 < p\le\infty$, then we understand $\widehat{\omega}$ in a distributional sense of  $S'(\R)$.

For ƒ$\sigma>0$,  we denote by $\bp$  the Bernstein space of  all $F\in L^p(\R)$ such that $F$ is bandlimited to $[-\sigma, \sigma]$.  The space $\bp$ is equipped with the norm
\[
\|F\|_{p}=\Bigl(\int_{\R}|F(x)|^p\,dx\Bigr)^{1/p} \ \ {\text{for}} \ \  \ 1\le p<\infty  \qquad and \qquad
\|F\|_{\infty}={\text{ ess supp}}_{x\in\R} |F(x)|.
\]
  By the Paley-Wiener-Schwartz theorem (see \cite[p. 68]{6}), any   $F\in \bp$ is infinitely differentiable on $\R$ and  has an extension onto the complex plane $\Co$  to an entire function of exponential type at most ƒ$\sigma$.  Note that $1\le p \le r \le\infty$  implies (\cite[p. 49]{6}, Lemma 6.6)
  \[
\be\subset \bp\subset B_{\sigma}^r\subset B_{\sigma}^{\infty}.
\]

There are several sampling expansion formulas for functions from $\bp$ like Whittaker-Kotel'nikov-Shannon  sampling formulas. For example, if $ F \in\bp$, $1\le p\le 2$, then $F$ can be expanded as the following  formula with derivatives  (see \cite[p. 60]{5})
\begin{equation}
F(x)=\sum_{n\in \Z}\biggl(F\biggl(\frac{2\pi}{\sigma}n\biggr)+F'\biggl(\frac{2\pi}{\sigma}n\biggr)\biggr(x-\frac{2\pi}{\sigma}n\biggr)\biggl)
\Biggl(\frac{\sin\biggl(\frac{\sigma}{2}x -\pi n\biggr)}{\frac{\sigma}{2}x -\pi n}\Biggl)^2.
\end{equation}
This series converges absolutely and uniformly on any compact subset of $\Co$.   For  $F\in \be$,  $x\in\R$, and   $a>0$, the Poisson summation formula reads (see, e.g.,  \cite[p. 63]{5} and \cite[p. 509]{2})
 \begin{equation}
\sum_{n\in\Z}F(x+an)=\frac1{a}\sum_{k\in\Z}\widehat F\biggl(\frac{2\pi}{a}k\biggr)e^{i\frac{2\pi x}{a}k},
\end{equation}
where both sums converge absolutely.

{\it Proof of Proposition 1.2} \
Define $\omega=\widehat{\mu}-\widehat{\eta}$. Then $\omega\in B(\R)$ and $\omega$ is compactly supported.  Using  Wiener's local theorem that the local belonging to $B(\R)$ is equivalent to the local belonging to $A(\R)=\{\widehat{\varphi}: \varphi\in L^1(\R)\}$ (see \cite[p. 258]{15}), we have that  $\mu -\eta \in L^1(\R)$. Therefore,  $\mu_s=\eta_s$  and Proposition 1.2 is proved.

{\it Proof of Theorem 1.3.} \
By Bochner's theorem, there is a probability measure $\mu$ such that $g=\widehat{\mu}$. The function $f-g$ is  continuous on $\R$ and supported on $[-\sigma,\sigma]$.  Therefore,   $(f-g)\in L^1(\R)$ and
\begin{equation}
\varphi- \mu= \widecheck{(f-g)}\in C_0(\R),
\end{equation}
where $C_0(\R)$ is the usual space of continuous functions on $\R$ that vanish at infinity. Hence,  (2.3) implies that $\mu$ is absolutely continuous with respect to the Lebesgue measure. Let $\psi\in L^1(\R)$ denote the density of $\mu$. Define $\zeta=\varphi-\psi$. According to the Bochner-Wiener-Schwartz theorem,  we conclude that $\zeta\in\be$. Moreover,
\begin{equation}
\zeta\le \varphi\quad {\text{on}}\quad \R
\end{equation}
and
\begin{equation}
\int_{\R}\zeta(x)\,dx=0.
\end{equation}
We claim that
\begin{equation}
\int_{E}\zeta(x)\,dx=0.
\end{equation}
To verify the claim, first let us define  $E_n=an+E$, $n\in\Z$. Then (1.2) and  (2.4) imply that
\begin{equation}
\int_{E_n}\zeta(x)\,dx\le 0
\end{equation}
for all $n\in\Z$. For $F=\zeta$, using the Poisson summation formula (2.2), we get
\begin{equation}
a\sum_{n\in\Z} \zeta(x+an)=\sum_{m\in\Z}\widehat{ \zeta}\Bigl(\frac{2\pi m}{a}\Bigr)e^{i\frac{2\pi m}{a}x}.
\end{equation}
By combining  the condition (1.3) with the fact that the continuous function $\widehat{ \zeta}$ is supported on $[-\sigma,\sigma]$,  we have that $\widehat{ \zeta}(2\pi m/a)=0$ for $m\in \Z\setminus\{0\}$. Moreover, we conclude from (2.5) that  $\widehat{ \zeta}(0)=\int_{\R}\zeta(x)\,dx=0$. Therefore, (2.8) implies that
\begin{equation}
\sum_{n\in\Z} \zeta(x+an)=0
\end{equation}
for all $x\in\R$.  The series on the left-hand side of (2.2) and (2.9) converge in $L^1[0,a]$ (see \cite[p. 509]{2}). Since $E\subset [0, a)$, it follows from (2.9)  that
 \[
0=\int_{E}\Bigl(\sum_{n\in\Z} \zeta(x+an)\Bigr)\,dx=\sum_{n\in\Z}\Bigl(\int_{E} \zeta(x+an)\,dx\Bigr)=\sum_{n\in\Z}\Bigl(\int_{E_n} \zeta(x)\,dx\Bigr).
\]
Combining this with (2.7), we see that $\int_{E_n} \zeta(x)\,dx=0$ for all $n\in\Z$. This proves our claim (2.6), since $E=E_0$.

Now we claim  that
\begin{equation}
 \zeta(x)=0\quad {\text{ for \ all}}\quad x\in E.
\end{equation}
To that end, we decompose $E$ into a disjoint union of its parts $E^{+}=\{x\in E: \zeta(x)>0\}$,  $E^{-}=\{x\in E: \zeta(x)<0\}$, and $E^{0}=\{x\in E: \zeta(x)=0\}$. We need only check that $E^{+}=E^{-}=\emptyset$. Indeed, using (1.2) and (2.4), we get
\[
\int_{E^{+}}\zeta(x)\,dx\le \int_{E^{+}}\varphi(x)\,dx\le  \int_{E}\varphi(x)\,dx=0.
\]
 This proves that $E^{+}=\emptyset$, since $\zeta$ is continuous on $\R$ and therefore   $E^{+}$ is an open subset of $\R$. Now it follows directly from  (2.6) that $\int_{E^{-}}\zeta(x)\,dx=0$. We have therefore $E^{-}=\emptyset$, which proves our claim (2.10).

Finally, (2.10) shows that the entire function $\zeta$ vanishes on the set $E$ of positive Lebesgue measure. Thus, by the uniqueness theorem for analytic functions, we have that $\zeta$ is the zero function or $\varphi\equiv \psi$. Theorem 1.3  is proved.

\begin{remark}
Suppose that  $f : \R\to\Co$ is the characteristic function of a   density $\varphi$. The proof of Theorem 1.3 implies that  in order   to exist  a  characteristic function $g$
  such that $f=g$ on $U_{\sigma}$ but $f\not\equiv g$, it is necessary and sufficient that there is    $\zeta\in \be$, $\zeta\not\equiv 0$, satisfying (2.4) and (2.5).
\end{remark}

{\it Proof of Corollary 1.4} \
Choose any $a\in \R$ such that
\begin{equation}
 |S_{\varphi}|< a\le \frac{2\pi}{\sigma}.
\end{equation}
Let $h_a: \R\to \R/a\Z$ denote the natural homomorphism of $\R$ onto the quotient group $ \R/a\Z$. We can identify  $ \R/a\Z$    with $[0,a)$  in the usual way. Then let us define
\begin{equation}
 E=[0,a)\setminus h_a(S_{\varphi}).
\end{equation}
We claim that $E$ satisfies the hypotheses of Theorem 1.3.  Indeed, $E$ is a measurable subset of $[0,a)$, since $ S_{\varphi}$ is measurable and $h_a$ is an open map. Next, by the definition of $h_a$  and  (2.11, we have
\[
| h_a(S_{\varphi})|\le |S_{\varphi}|<a.
\]
Hence $|E|>0$.  Next,  suppose that there exists  a  $x\in S_{\varphi}$ such that $x\in E+a\Z$. Then $h_a(x)\in h_a(S_{\varphi})$ and together $h_a(x)\in h_a(E+a\Z)=E$. This contradicts the definition of $E$ in (2.12). Therefore, we see that $( E+a\Z)\cap S_{\varphi}=\emptyset$. This implies (1.2), proving  our claim. Finally, using Theorem 1.3, we complete the proof of Corollary 1.4.

Before proving Theorem 1.5, we shall need the following lemma.
\begin{lemma}.
Let $a, b\in (-\infty, \infty)$, $a<b$.  Suppose that  $F\in\be$, $F\not\equiv 0$, satisfies
\begin{equation}
 \int_{\R} F(x)\,dx=0
\end{equation}
and
\begin{equation}
 F(x)\le \chi_{[a,b]}(x)
\end{equation}
for all $x\in\R$, where $\chi_{[a,b]}$ is the indicator function of $[a,b]$. If $b-a=2\pi/\sigma$, then there exists a $\tau\in(0, 2\pi/\sigma]$  such that
\begin{equation}
 F(x)=\tau\biggl[\frac1{x-a}+\frac{\sigma}{2\pi-\sigma(x-a)}\biggr]\sin^2\frac{\sigma(x-a)}{2},
 \end{equation}
 $x\in\R$.
\end{lemma}

{\it Proof.} \
Using the change of variables $x\to x+a$, we may  assume  further without loss of generality, that $a=0$. Then  $b=2\pi/\sigma$. Now  we claim that
\begin{equation}
 F\biggl(\frac{2\pi}{\sigma}n \biggr)=0
 \end{equation}
 for all $n\in\Z$. Indeed, (2.14) implies clearly  that
  \begin{equation}
   F\biggl(\frac{2\pi}{\sigma}n \biggr)\le 0, \quad n\in\Z\setminus\{0,1\}.
    \end{equation}
   Moreover, since   $F$ is continuous on $\R$, we see   that $F$ satisfies (2.17)   also at the endpoints $0$ and $2\pi/\sigma$ of $[0, 2\pi/\sigma]$. Combining this with (2.13) and using the following summation formula for  $F\in \be$ (see \cite[p.  122]{3})
   \[
   \sigma\int_{\R}F(x)\,dx=\sum_{n\in\Z}F\biggl(\frac{2\pi}{\sigma}n\biggr),
\]
we obtain our claim (2.16).

By (2.14),  we see that $F(x)\le 0$ on $\R\setminus [0, 2\pi/\sigma]$. Then (2.16) implies that $F$ attains a maximum (local) at every point $ 2n\pi/\sigma$,   $n\in\Z\setminus\{0,1\}$. Hence
 \begin{equation}
   F'\biggl(\frac{2\pi}{\sigma}n \biggr)=0
    \end{equation}
for $n\in\Z\setminus\{0,1\}$. Now substituting (2.16) and (2.18) into (2.1), we have
\begin{equation}
   F(x)= xF'(0)\biggl(\frac{\sin\frac{\sigma x}{2}}{\frac{\sigma x}{2}}\biggr)^2+\Bigl(x-\frac{2\pi}{\sigma}\Bigr)F'\Bigl(\frac{2\pi}{\sigma}\Bigr)\biggl(\frac{\sin\frac{\sigma x}{2}}{\frac{\sigma x}{2}-\pi}\biggr)^2.
\end{equation}
It is straightforward to check that this function belongs to $\be$  if and only if $F'(0)=-F'\Bigl(\frac{2\pi}{\sigma}\Bigr)$. Therefore, the function (2.19) reduces to
\[
  F(x)= \tau \biggl[\frac{1}{x}+\frac{\sigma}{2\pi -\sigma x}\biggr]\sin^2\frac{\sigma  x}{2}
  \]
with certain $\tau\neq 0$. Moreover, (2.14) implies that  $\tau> 0$. Then
\[
\max_{x\in [0, \ 2\pi/\sigma]} F(x)=F\Bigl(\frac{\pi}{\sigma}\Bigr)=\frac{2\sigma \tau}{\pi}.
\]
Thus,  it follows from (2.14) that $\tau$ can be  an arbitrary  number such that $0<\tau\le \pi/(2\sigma)$. Lemma 2.2  is proved.

{\it Proof of Theorem 1.5.} \
Let
\[
\varepsilon=\frac12(\beta-\alpha)-\frac{\pi}{\sigma}.
\]
It follows   from (1.6)  that $\varepsilon>0$.  Clearly,  the length of  $[\alpha+\varepsilon, \beta-\varepsilon]$ is $2\pi/\sigma$.  Define
\[
\varrho=\min\{\varphi(x): \ x\in [\alpha+\varepsilon, \beta-\varepsilon]\}.
\]
The quantity $\varrho$  is well defined, since now $\varphi$ is continuous. Moreover, (1.5) gives that $\varrho>0$.  Assume  that  $F$ is the function (2.15) with $a=\alpha+\varepsilon$  and an  arbitrary fixed $\tau\in(0, 2\pi/\sigma]$. Set $G=\varrho F$.  Then $G\in\be$ and  $G\not\equiv 0$. In addition,
\[
G(x)\le \varphi(x), \quad x\in\R,
\]
since the function (2.15) satisfies $F(x)\le 0$ for all $x\in \R\setminus [a, a+2\pi/\sigma]$. On the other hand, (2.13) implies that $\int_{\R}G=\varrho\int_{\R}F=0$. Therefore, we see that $\varphi-G$ is a probability density. Set $ g=\widehat{\varphi-G}$. Then $g$ is the characteristic function such that  $f(x)=\widehat{\varphi}(x)=\widehat{(\varphi-G)}(x)=g(x)$ for
$x\in U_{\sigma}$ but $f\not\equiv g$,  since $G\not\equiv 0$.  Theorem 1.5  is proved.

{\it Proof of Proposition 1.8.}

To simplify the proof,  we will show that there exist continuously differentiable functions $\varphi\in C'(\R)$  for which  the requirements (i) or (ii) of Proposition 1.8  are satisfied.  Therefore,  throughout this proof we assume that $\varphi\in C'(\R)$.

By  Remark 2.1, we see that  each characteristic function $g:\R\to\Co$ such that $g=\widehat{\varphi}$ on $U_{\sigma}$ has the form $\widehat{\varphi}+\widehat{\zeta}$, where $\zeta$  is any function in $\be$ satisfying (2.4) and (2.5). Note that now we do not  require that $g\not\equiv f$, i.e., it can be that  $\zeta\equiv 0$. For any such    $\zeta$ and
\[
\varrho=\max\Bigl\{\varphi(x): x\in \Bigl[\alpha, \alpha+\frac{2\pi}{\sigma}\Bigr]\Bigr\},
\]
let us define $F=\zeta/\varrho$. Then  $F$  satisfy all the conditions of  Lemma 2.2  except the hypothesis that $F\not\equiv 0$. Therefore, we may assume that $F$ is the function (2.15) with $a=\alpha$ and certain $\tau\in (0, 2\pi/\sigma]$ or $\tau=0$ (the case $\zeta\equiv 0$).  Then a straightforward calculation shows that $F(\alpha)=F(\alpha+2\pi/\sigma)=0$ and
\begin{equation}
F'(\alpha)=-F'\Bigl(\alpha+\frac{2\pi}{\sigma}\Bigr)= \frac{\tau\sigma^2}{4}.
\end{equation}
(i) \ Let $\varphi\in C'(\R)$ be a density function such that  $S_{\varphi}=(\alpha, \alpha+2\pi/\sigma)$.  Suppose that  $\varphi$ satisfies  at least one of the following conditions
\begin{equation}
\varphi'(\alpha)= 0\quad{\text{or}}\quad \varphi'\Bigl(\alpha+\frac{2\pi}{\sigma}\Bigr)= 0.
\end{equation}
Then (2.4) implies that $\zeta'(\alpha)\le 0$ or $\zeta'(\alpha+2\pi/\sigma)\ge 0$, since $\zeta(\alpha)=F(\alpha)=0$ and $\zeta(\alpha+2\pi/\sigma))=F(\alpha+2\pi/\sigma)=0$. Therefore, by the definition  of $F$, we obtain from  (2.20)  that $\tau=0$,  i.e.,   $\zeta\equiv 0$. Hence, any such $\varphi$ satisfying  at least one of the conditions (2.21) also satisfies the statement (i) of Proposition 1.8.

(ii) \ Suppose that $\varphi\in C'(\R)$ is any density function   that satisfies  $S_{\varphi}=(\alpha, \alpha+2\pi/\sigma)$  and
\[
\varphi'(\alpha)\neq 0\quad {\text{ and}}\quad \varphi'\Bigl(\alpha+\frac{2\pi}{\sigma}\Bigr)\neq  0.
\]
Now we claim that there  exists $\zeta\in\be$, $\zeta\not\equiv 0$,   such that  $\zeta$ satisfies (2.4) and (2.5). Indeed, it is easily see that there exist $\varepsilon>0$ and $\tau\in(0, 2\pi/\sigma]$  for which the function (2.15) with $a=\alpha$ satisfies
\[
F(x)\le \varphi(x)\quad {\text{for}}\quad x\in \Bigl[\alpha,\alpha+\varepsilon\Bigr]\cup \Bigl[\alpha+\frac{2\pi}{\sigma}-\varepsilon, \alpha+\frac{2\pi}{\sigma}\Bigr].
\]
Since  $S_{\varphi}=(\alpha, \alpha+2\pi/\sigma)$ and $\varphi\in C'(\R)$, we may suppose (passing to a smaller positive value of $\tau$ in (2.15) if necessary) that the condition $F(x)\le \varphi(x)$  is satisfied for all $x\in\R$. Finally, if we  set $\zeta=F$ and $g=\widehat{\varphi}+\widehat{\zeta}$, then these functions satisfy the statement (ii) of our proposition.

}}
\end{document}